\documentclass[12pt]{article}

\usepackage[latin2]{inputenc}
\usepackage{amsmath}
\usepackage{amssymb}
\usepackage{amsfonts}
\usepackage{stmaryrd}
\usepackage{theorem}
\usepackage{graphicx}
\usepackage{epstopdf}
\usepackage[parfill]{parskip}    
\usepackage{lineno}
\usepackage{float}
\usepackage{hyperref}
\hypersetup{pdftex,colorlinks=true,allcolors=blue}
\usepackage{hypcap}
\usepackage[super]{natbib}
\usepackage{authblk}

\newtheorem{thm}{Theorem}[section]
\newtheorem{lem}{Lemma}[section]
\newtheorem{cor}{Corollary}[section]

\newtheorem{dfn}{Definition}[section]

\newtheorem{alg}{Algorithm}[section]

\newcommand{\keywords}[1]{\par\noindent{\small{\em Keywords\/}:\ \ #1}}
\newcommand{\mscclass}[1]{\par\noindent{\small{\em MSC class\/}:\ \ #1}}

\newcommand*{\mH}{\mathrm{H}}

\newcommand*{\adr}{\mathrm{adr}}
\newcommand*{\diam}{\mathrm{diam}}
\newcommand*{\mymod}{\ \mathrm{mod}\ }
\newcommand*{\Conv}{\mathrm{Conv}}

\newcommand*{\cf}{\ensuremath{\varphi}}

\newcommand*{\ve}{\ensuremath{\varepsilon}}
\newcommand*{\vt}{\ensuremath{\vartheta}}
\newcommand*{\vl}{\ensuremath{\lambda}}

\newcommand*{\vm}{\ensuremath{\mu}}
\newcommand*{\N}{\ensuremath{\mathbb{N}}}

\newcommand*{\Q}{\ensuremath{\mathbb{Q}}}
\newcommand*{\R}{\ensuremath{\mathbb{R}}}
\newcommand*{\C}{\ensuremath{\mathbb{C}}}
\newcommand*{\prf}{\textbf{Proof}\ \ }

\newcommand*{\sqr}{\ensuremath{\square}}

\newcommand*{\frar}{\ensuremath{\shortrightarrow}}

\newcommand*{\mcN}{\ensuremath{\mathcal{N}}}
\newcommand*{\mcT}{\ensuremath{\mathcal{T}}}
\newcommand*{\mcP}{\ensuremath{\mathcal{P}}}
\newcommand*{\mcA}{\ensuremath{\mathcal{A}}}
\newcommand*{\co}{\ensuremath{\circ}}
\newcommand*{\tti}{\ensuremath{\rightarrow\infty}}
\newcommand*{\mes}{\ensuremath{\varnothing}}

\begin{document}

\title{On Intersecting IFS Fractals with Lines}
\author{József Vass}
\affil{Department of Applied Mathematics\\ University of Waterloo\\ 200 University Avenue West\\ Waterloo, ON, N2L 3G1, Canada\\ jvass@uwaterloo.ca}
\date{}
\maketitle

\begin{abstract}
\noindent IFS fractals - the attractors of Iterated Function Systems - have motivated plenty of research to date, partly due to their simplicity and applicability in various fields, such as the modeling of plants in computer graphics, and the design of fractal antennas. The statement and resolution of the Fractal-Line Intersection Problem is imperative for a more efficient treatment of certain applications. This paper intends to take further steps towards this resolution, building on the literature. For the broad class of hyperdense fractals, a verifiable condition guaranteeing intersection with any line passing through the convex hull of a planar IFS fractal is shown, in general $\R^d$ for hyperplanes. The condition also implies a constructive algorithm for finding the points of intersection. Under certain conditions, an infinite number of approximate intersections are guaranteed, if there is at least one. Quantification of the intersection is done via an explicit formula for the invariant measure of IFS.\footnote{The first draft of the paper was shared on Dec. 23, 2011. The second draft was submitted on Dec. 25, 2012 and was accepted for publication on Jun. 29, 2014 in the journal Fractals \copyright\ 2014 World Scientific Publishing Company \url{http://www.worldscientific.com/worldscinet/fractals}.}
\ \\
\mscclass{28A80 (primary); 37F99, 52A35 (secondary).}
\end{abstract}

\newpage
\section{Introduction} \label{s00} 
Falconer\cite{bb00003} surveys the properties of fractals under projections. When the projection is carried out in $\R^2$ or $\R^3$ onto a line or plane, one might consider the resulting set as the ``shadow'' of the fractal, which may be analyzed for its own dimension. Furthermore, we might ask ``how many'' fractal points are projected to a certain point on the line, or if there are any at all, in essence inquiring about the distribution of such a projection. So the projection problem breaks down into two main questions, since the directional ray of projection can be thought of as an intersecting line:\\
(1) The Fractal-Line Intersection Problem: Given a line and an IFS fractal in the plane, do they intersect?\\
(2) If they intersect, how many points of intersection are there? How about within some $\ve>0$ neighbourhood of the line?\\
This paper intends to resolve these general questions, while hinting at their relevance for applications. For certain broad classes of IFS fractals - hyperdense, or specifically chain fractals - the shadow is shown to be always filled in, no matter where the light shines from, proven in general in $\R^d$. Furthermore, we show that the segment shadow in $\R^2$ receives an infinite number of projected points, in any $\ve>0$ subinterval of the segment. These properties may make some of these potentially disconnected fractals ideal for 2D fractal antenna design, or as light-absorbing tree crowns (which could be considered 3D fractal antennas).

Recently Mendivil and Taylor\cite{ba00018} approached these questions from a projectional point of view. Defining a certain parametrized class of planar IFS fractals, they wish to guarantee that the shadow in all directions is a segment. They prove that this holds for some domain of parameters. In other words, for this specific class of planar fractals, any line or ray of light that intersects the convex hull, also intersects the fractal, thereby contributing to its shadow. We examine the problem further for the broadest possible class of attractors in $\R^d$ called hyperdense fractals, and we also introduce the verifiable subclass of chain fractals.

The Fractal-Line Intersection Problem is relevant to a number of applications, among which we mention ray tracing in computer science, the design of fractal antennas in engineering, and the study of tree crown density for light absorption in botany and forestry. In computer science, ray tracing involves the shading of an object in virtual space, which is detailed in Hart and DeFanti\cite{bc00005} in regards to 3D IFS fractals, as well as numerous other publications. Fractal antennas are flat metal antennas with an IFS fractal layout that must be optimized for the amount of material used versus the efficiency of signal reception. These antennas were introduced by Cohen et al.\cite{ba00030, ba00019}. The study of plant growth and tree crowns for light absorption are vast fields, for which see Prusinkiewicz and Lindenmayer\cite{bb00005} and Zeide\cite{ba00020}. Last but not least, the pioneering inspirational work of Mandelbrot\cite{bb00001} must be emphasized.

Further research into the projection of fractals has been carried out by Besicovitch\cite{ba00021} and Federer\cite{ba00022} examining s-sets; Marstrand\cite{ba00024}, Kaufman\cite{ba00025}, and Mattila\cite{ba00026, bb00006} showing projection theorems for arbirary sets in $\R^d$; and Davies\cite{ba00023}, Falconer et al.\cite{ba00027, ba00028}, and Howroyd\cite{ba00029} giving results for box and packing dimensions. These efforts are all summarized in the expository book by Falconer\cite{bb00003}, which also provides an introduction to IFS fractals.

\section{Preliminary Concepts} \label{s01}

The attractors of Iterated Function Systems - IFS fractals - were pioneered by Hutchinson\cite{ba00007}, further discussed by Barnsley and Demko\cite{ba00017}, and may be the most elementary fractals possible. They are the attractors of a finite set of affine linear contraction mappings on $\R^d$ - the ``function system'' - which when combined and iterated to infinity, converges to an attracting limit set, the IFS fractal itself. We begin by defining the generating maps of an IFS fractal, and go on to stating its existence and uniqueness.

\subsection{IFS Fractals} \label{s0101}

\begin{dfn} \label{s010101}
Let an affine contraction mapping (briefly: contraction) $T:\R^d\shortrightarrow\R^d$ be defined for all $z\in\R^d$ as $T(z):=p+M(z-p)$ where $p\in\R^d$ is the fixed point of $T$, and the invertible $M\in\R^{d\times d}$ is the factor of $T$, with $\|M\|_2<1$.
\end{dfn}

When $T$ is a similitude over the complex plane ($d=2$), it may be written in the form $T(z):=p+\varphi(z-p)$ where $z, p\in\C$ and $\varphi=\lambda e^{\theta i}\in\C$, with $\lambda\in(0,1)$ and $\theta\in(-\pi,\pi]$. In higher dimensions, similitudes are $M=\vl R$ with a unitary $R$.

\begin{dfn} \label{s010102}
Let an affine contractive $n$-map iterated function system (briefly: IFS) be defined as a finite set of contractions, and denoted as $\mathcal{T}:=\{T_1,\ldots,T_n\},\ n\in\N$. Further denote $\mcN:=\{1,\ldots,n\},\ \mcP:=\{p_1,\ldots,p_n\},\ \Phi:=\{M_1,\ldots,M_n\}$.
\end{dfn}

\begin{dfn} \label{s010103}
Let $\mathcal{T}=\{T_1,\ldots,T_n\}, n\in\N$ be an IFS. Define the Hutchinson operator $\mH$ belonging to $\mathcal{T}$ as
\[ \mH(S):=\bigcup_{k=1}^n T_k(S),\ \ T_k(S):=\{T_k(z): z\in S\},\ \mathrm{for\ any}\ S\subset\R^d \]
and call $\mH(S)$ the Hutchinson of the set $S$.
\end{dfn}

\begin{thm} \label{s010104}
For any IFS with Hutchinson operator $\mH$, there exists a unique compact set $F\subset\R^d$ such that $\mH(F)=F$. Furthermore, for any compact $S_0\subset\R^d$, the recursive iteration $S_{n+1}:=\mH(S_n)$ converges to $F$ in the Hausdorff metric.
\end{thm}
\prf
The proof follows from the Banach Fixed Point Theorem, once we show that $\mH$ is contractive in the Hausdorff metric over compact sets\cite{ba00007}. \sqr

\begin{dfn} \label{s010105}
Let the set $F$ in the above theorem be called a fractal generated by an IFS with Hutchinson operator $\mH$ (briefly: IFS fractal). Denote $\langle T_1,\ldots,T_n\rangle = \langle\mathcal{T}\rangle :=F$.
\end{dfn}

\subsection{The Address Set} \label{s0102}

The address set results from the iteration of the Hutchinson operator, and it is a way to label the location of each fractal point. Since we can start the iteration towards $F$ with any compact set, we often choose the primary fixed point, which is any point in $\mcP$ of our preference.

\begin{dfn} \label{s010201}
Let $\mcN^j:=\mcN\times\ldots\times\mcN$ be the index set to the $j$-th Cartesian power, and call this $j$ the iteration level. Then define the address set as $\mcA:=\{0\}\cup\bigcup_{j=1}^\infty\mcN^j$. For any $a\in\mcA$ denote its $k$-th coordinate as $a(k),\ k\in\N$. Let its dimension or length be denoted as $|a|\in\N$ so that $a\in\mcN^{|a|}$ and let $|0|:=0$. Define the map with address $a\in\mcA$ acting on any $z\in\R^d$ as the function composition $T_a(z):=T_{a(1)}\circ\ldots\circ T_{a(|a|)}(z)$. Let the identity map be $T_0:=Id$.
\end{dfn}

The above definition of the address set merges the two standard definitions in use. Common ways include identifying addresses with decimal number representations with a certain basis, and strings of letters usually for a small number of maps. This formal language representation accounts for finite addresses by using vectors of numbers, preferable to the string formalism when $|\mcN|>24$.

\begin{thm} \label{s010202}
For any primary fixed point $p\in\mcP$ we have
\[ F=\lim_{N\rightarrow\infty} \mH^N(\{p\})=\mathrm{Cl}\{T_a(p):a\in\mcA\}=\mathrm{Cl}\{T_a(p_k):a\in\mcA,\ p_k\in\mcP\} \]
and we call this the address generation of $F$.
\end{thm}
\prf
The proof follows from Theorem \ref{s010104} with either of the initial sets $\{p\}$ or $\mcP$. \sqr

\begin{dfn} \label{s010203}
Let the address $\adr_p(f)$ of a fractal point $f\in F$ with respect to a primary fixed point $p\in\mcP$ be the shortest address $a\in\mcA$ for which $T_{a}(p)=f$ (if two such addresses exist equal in length, then take the lexicographically lower one).
\end{dfn}

\subsection{Hyperdense Fractals} \label{s0103}

\begin{dfn} \label{s010301}
Let an IFS fractal in $\R^d$ be hyperdense if any hyperplane that intersects its convex hull, also intersects the Hutchinson of its convex hull. Let an IFS fractal be a chain fractal, if the Hutchinson of its convex hull is connected.
\end{dfn}

\begin{thm} \label{s010302}
Chain fractals are hyperdense.
\end{thm}
\prf
Denote the fractal $F=\langle T_1,\ldots,T_n\rangle$ and its convex hull $C_F:=\Conv(F)$. First we see that for any $T\in\mcT$ we have $T(\Conv(S))=\Conv(T(S))$ since $T$ is affine.\\
Note also that for any $S_1,\ldots,S_N\subset\R^d$
\[ \Conv\left(\bigcup_{k=1}^N S_k\right) = \Conv\left(\bigcup_{k=1}^N \Conv(S_k)\right) \]
since the convex combination of convex combinations, is a convex combination. Considering the fact that $\mH(F)=F$ as well as the above ideas, we have that
\[ \Conv(F)=\Conv\left(\bigcup_{k=1}^n \Conv(T_k(F))\right)=\Conv\left(\bigcup_{k=1}^n T_k(\Conv(F))\right) \]
meaning that $C_F=\Conv(\mH(C_F))$. Since $F$ is compact so is $C_F$ and thus $\mH(C_F)$, since $\mH$ is continuous.\\
We now turn to showing that $F$ is hyperdense. Let us take any hyperplane $L\subset\R^d$ that intersects $C_F=\Conv(\mH(C_F))$ in some point $q=\vm h_1+(1-\vm)h_2,\ h_{1,2}\in \mH(C_F),\ \vm\in[0,1]$. If $\vm\in\{0,1\}$ then $q\in \mH(C_F)$ so we are done. Otherwise $L$ separates the space into two half spaces, with $h_1$ in one and $h_2$ in the other. Since $F$ is a chain fractal, we know that $\mH(C_F)$ is connected, thus it is path-connected, so there is a path $\gamma\subset \mH(C_F)$ connecting $h_1$ and $h_2$. Since $h_{1,2}$ are on separate sides of $L$, we must have that $L\cap\gamma\neq\mes$. This can be shown by parametrizing $\gamma:[0,1]\frar\R^d$ and writing the plane as $L=\{z\in\R^d: \langle a,z\rangle= b\}$ for some $a\in\R^d,\ b\in\R$. Denoting $f(t):=\langle a,\gamma(t)\rangle -b$ we have that $f(0)f(1)<0$, so by Bolzano's theorem $f$ must have a root $t_0\in (0,1)$, implying that $\gamma(t_0)\in L\cap\gamma\subset L\cap \mH(C_F)$. Therefore by $L\cap\gamma\neq\mes$ and $\gamma\subset \mH(C_F)$ we have that $L\cap \mH(C_F)\neq\mes$. \sqr

\begin{lem} \label{s010303}
For a hyperdense fractal $F=\langle T_1,\ldots,T_n\rangle$ and any address $a\in\mcA$, if a hyperplane $L$ intersects $T_a(C_F)$ then it also intersects $T_a(\mH(C_F))$.
\end{lem}
\prf
It is clear that since the $M_k\in\R^{d\times d}$ factors in $T_k$ are invertible, the inverses of the maps are $T_k^{-1}(z)=p_k+M_k^{-1}(z-p_k)$. Thus $T_a^{-1}$ also exists, and it is also an affine mapping, so it takes the hyperplane $L$ into another hyperplane $L'$. Thus $L\cap T_a(C_F)\neq\mes$ is equivalent to $T_a^{-1}(L)\cap C_F\neq\mes$, which by the hyperdensity of $F$ implies that $T_a^{-1}(L)\cap \mH(C_F)\neq\mes$, and so $L\cap T_a(\mH(C_F))\neq\mes$. \sqr

\section{Fractal-Line Intersection} \label{s02}

\subsection{Exact Intersection} \label{s0201}

\begin{thm} \label{s020101}
A hyperplane intersects a hyperdense fractal if and only if it intersects its convex hull. This equivalence holds only if the fractal is hyperdense.
\end{thm}
\prf
The proof is based on Cantor's Intersection Theorem and the address generation of $F=\langle T_1,\ldots,T_n\rangle\subset\R^d$ in Theorem \ref{s010202}. We show that the hyperdensity of $F$ implies a decreasing sequence of compact sets, which tend to a point by Cantor's Intersection Theorem. The index sequence itself will correspond to an address, which in the limit locates a fractal point, since the fractal is the closure of all possible addresses. Let us now begin the proof.\\
If a hyperplane $L$ intersects $F$, it must clearly intersect $C_F:=\Conv(F)$. On the other hand, if $L$ intersects $C_F$, by the fractal's hyperdensity, it also intersects $\mH(C_F)$. So $L$ must intersect $T_{k_1}(C_F)$ for some $k_1\in\mcN$. Let this intersection be denoted as $I_1:=L\cap T_{k_1}(C_F)\subset\R^d$. Then $I_1$ is compact, because $C_F$ is compact.\\
Since $T_{k_1}(\mH(C_F))=\bigcup_{k=1}^n T_{k_1}\co T_k(C_F)$, according to Lemma \ref{s010303} with $a=(k_1)$, the fact that $L$ intersects $T_{k_1}(C_F)$ implies that it also intersects $T_{k_1}(\mH(C_F))$ and thus $T_{k_1}\co T_{k_2}(C_F)$ for some $k_2\in\mcN$. Once again denoting $I_2:=L\cap T_{k_1}\co T_{k_2}(C_F)$ we have that this set is compact, and $I_2\subset L\cap T_{k_1}(C_F)=I_1$.\\
Continuing to apply Lemma \ref{s010303} in the above recursive procedure, by induction we get a strictly monotonically decreasing sequence (since $T_k,\ k\in\mcN$ are contractive) of compact sets in $\R_d$: $I_1\supset I_2\supset\ldots\supset I_j\supset\ldots$ each with a corresponding address composition: $T_{k_1}\co\ldots\co T_{k_j}$. According to Cantor's Intersection Theorem $\bigcap_{j=1}^\infty I_j\neq\mes$ and it contains a single point $f\in\R^d$, since the address composition contracts to a point in the limit.\\
Starting with any $p\in \mcP\cap C_F$ we have that $f=\lim_{j\rightarrow\infty} T_{k_1}\co\ldots\co T_{k_j}(p)\in F$ by the address generation of $F$, so $L$ intersects $F$ in $f$. Note that such an intersection may not be unique, since in our recursive proof, we only chose one indexed set in the Hutchinson union at each step, though $L$ may intersect multiple.\\
The above is under the condition that $F$ is hyperdense. If it is not, then by definition there is a hyperplane which intersects $C_F$ but not $\mH(C_F)$. Since $F\subset \mH(C_F)$ this hyperplane will not intersect $F$, countering the equivalence. \sqr

\begin{cor} \label{s020102}
A line intersects a chain fractal in $\R^2$ iff it intersects its convex hull.
\end{cor}

\begin{alg} \label{s020103}
(Fractal-Line Intersection) Let $F=\langle T_1,\ldots, T_n\rangle$ be a hyperdense fractal (possibly a chain fractal), and assume that its convex hull $C_F$ is known explicitly. Furthermore, let $L$ be a line, and $\ve>0$ an arbitrary stopping parameter.\vspace{-0.5em}\begin{tabbing}
Step 0: $a:=(0),\ I:=\mes$\\
Step 1: \=If \=$L\cap T_a(C_F)\neq\mes$ then\\
\> \> For\=\ $k=1$ to $n$\\
\> \> \> If \=$L\cap T_a\co T_k(C_F)\neq\mes$ then\\
\> \> \> \> $a:=(a,k)$\\
\> \> \> \> If \=$\diam(L\cap T_a(C_F))\ge\ve$ then\\
\> \> \> \> \> GoTo Step 1\\
\> \> \> \> Else\\
\> \> \> \> \> GoTo Step 2\\
\> \> \> \> End\\
\> \> \> End\\
\> \> End\\
\> End\\
Step 2: $I:=I\cup (L\cap T_a(C_F))$
\end{tabbing}
\end{alg}
\prf
It follows directly from the proof of Theorem \ref{s020101}. The algorithm finds a union of intervals - each with length below $\ve$ - that covers the points of intersection of the line and the fractal. \sqr

\subsection{Approximate Intersection} \label{s0202}

We proceed to the second question posed earlier, the number of intersections within some $\ve>0$ accuracy. These results are also relevant for applications, such as the signal-reception efficiency of fractal antennas. If an antenna is designed to be hyperdense - possibly a chain fractal - then not only will it intercept all signal planes crossing its convex hull - making it space-efficient - but will do so an infinite number of times, in an approximate sense.

\begin{thm} \label{s020201}
If an open set has a common point with an IFS fractal, then it has an uncountably infinite number.
\end{thm}
\prf
Let $S\subset\C$ the open set and $f\in F\cap S$. We may suppose that $f$ has a finite address $a\in\mcA,\ |a|<\infty$, since in the address generation of $F$, the fractal points with a finite address are dense in $F$, which is their closure according to Theorem \ref{s010202}. So we may replace $f$ with another fractal point $f'$ that has a finite address, and which is close enough to $f$ to be still be an element of $S$. Thus we may suppose that $f$ has a finite address, with a corresponding map $T_a$ and primary fixed point $p$, meaning $f=T_a(p),\ p\in\mcP$. Denote the contraction belonging to $p$ as $T\in\mcT$.\\
Let $\ve'>0$ be the radius of some ball centered at $f$ - denoted by $B':=B(f,\ve')$ - that is still contained in $S$. If we show that there are an uncountably infinite number of fractal points in $B'$, then that implies the theorem.\\
Let us map back $f=T_a(p)$ and $B'$ by $T_a^{-1}$ to $p=T_a^{-1}(f)$ and $B'':=T_a^{-1}(B')=B(p,\ve'')$ respectively. Here $\ve''=\frac{\ve'}{|\cf_a|}$ where $\cf_a$ is the product of the factors of the contractions in $\mcT$ making up $T_a$. Then since $F$ is compact, we may map it iteratively by $T$ until it is contained in $B''$, that is $\exists k\in\N: T^k(F)\subset B''$, and clearly $T^k(F)$ has an uncountably infinite number of points.\\
Mapping it all back by $T_a$, we have that $T_a\co T^k(F)\subset F\cap B'\subset F\cap S$. Thus we have shown an uncountably infinite number of common points of $F$ and $S$. \sqr

\begin{cor} \label{s020202}
Suppose that some open $\ve>0$ translational neighborhood $L_{\ve}^t$ of a line $L$ contains a point in an IFS fractal $F$, meaning
\[ \exists \ve>0, f\in F: f\in L_{\ve}^t:=\{z\in\C: d(L,z)<\ve\} \]
Then there are an uncountably infinite number of fractal points in $L_{\ve}^t$.
\end{cor}

\begin{cor} \label{s020203}
Suppose that some open $\ve>0$ angular neighborhood $L_{\ve}^a$ of a line $L$ around $q\in L$ contains a point in an IFS fractal $F$, meaning
\[ \exists \ve>0, f\in F: f\in L_{\ve}^a:=\{z\in\C: \angle(L,z-q)<\ve\} \]
Then there are an uncountably infinite number of fractal points in $L_{\ve}^a$.
\end{cor}

Both of the above corollaries follow directly from the theorem. The next theorem is of a slightly different nature, and is not implied by it. We restrict ourselves to IFS fractals in the complex plane, generated by an IFS of similitudes.

\begin{thm} \label{s020204}
Let $\mcT$ be an IFS with at least one map $T\in\mcT$ having a rotation angle $\vt$ for which $\frac{\vt}{2\pi}$ is irrational. Then any line $L$ that intersects $F=\langle\mcT\rangle$ in some point $f$ with a finite address, intersects it in at least a countably infinite number of points, with any $\ve>0$ angular accuracy around $f$, meaning
\[ \forall \ve>0\ \exists (f_k)_{k=1}^{\infty}\subset F: \angle(L,f_k-f)<\ve,\ k\in\N \]
\end{thm}
\prf
The address of $f$ is finite, meaning that $f=T_a(p),\ p\in\mcP,\ a\in\mcA,\ |a|<\infty$. Let us transform back $f$ and $L$ to $p\in\mcP$ by $T_a^{-1}$. Then we have that $p=T_a^{-1}(f)$, and denote $L':=T_a^{-1}(L)$ which is also a line. The $\ve$ angular neighborhood of $L$ around $f$ is transformed by $T_a^{-1}$ to an $\ve$ angular neighborhood of $L'$ around $p$. Clearly $p\in L'$ since $f=T_a(p)\in L=T_a(L')$. For the theorem to hold, it is sufficient to find an infinite number of fractal points within the angular neighborhood of $L'$, since we can map these points with $T_a$ to the angular neighborhood of $L$.\\
Let us choose any $q\in F,\ r\in L'$, with $\alpha:=\arg(q-p),\ \beta:=\arg(r-p)$. Then $\arg (T^k(q)-p)=(k\vt +\alpha)\mymod 2\pi$ and the iterations of $q$ by $T$ will be along a logarithmic spiral around $p$. We show that the iterations visit the $\ve$ angular neighborhood of $L'$ infinitely often.
\begin{figure}[H]
\centering
\includegraphics[width=200pt]{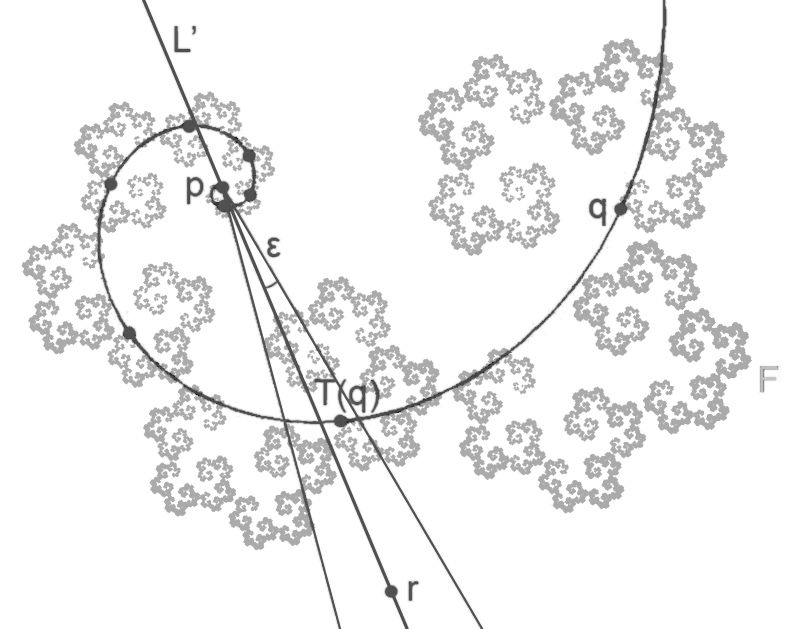}
\caption{Illustration of the angular neighborhood of $L'$ around $p$.}
\label{s020205}
\end{figure}
We have supposed that $\frac{\vt}{2\pi}$ is irrational, so by the Equidistribution Theorem, the sequence $(k\vt\mymod 2\pi)_{k\in\N}$ is uniformly distributed in $[0,2\pi)$, and thus it is also dense in this interval. So approximating the angle between $q$ and $r$ with respect to $p$, ie. $\beta-\alpha$, we have the following
\[ \forall j\in\N\ \exists k_j\in\N: ((k_j\vt +\alpha)-\beta)\mymod 2\pi < \frac{\ve}{j}<\ve \]
Therefore the sequence $(T^{k_j}(q))_{j\in\N}$ will be within the $\ve$ angular neighborhood of $L'$ with respect to $p$. Thus mapping the sequence back by $T_a$, it will be in the required $\ve$ angular neighborhood of $L$ with respect to $f$, and clearly $(T_a\co T^{k_j}(q))_{j\in\N}\subset F$. So with $f_j:=T_a\co T^{k_j}(q),\ j\in\N$ we have shown a countably infinite number of fractal points in the $\ve$ angular neighborhood of the intersecting line. \sqr

\begin{cor} \label{s020206}
Suppose a line intersects a planar IFS fractal in a point with a finite address, and has some $\ve>0$ angular neighborhood around the point that contains only a finite number of fractal points. Then all contractions in the IFS have a rotation angle $\vt$ for which $\frac{\vt}{2\pi}\in\Q$.
\end{cor}

The above corollary hints at the relevance of a certain class of planar IFS fractals, having roots of unity as rotation factors. They may hold a special place in the connectedness of IFS fractals. Examining such ``rational fractals'' further may prove to be a fruitful venture, and the case of primitive roots of unity may be even more worthwhile. Indeed these theorems seem to call for an investigation into the translational and angular distribution of IFS fractals.

\subsection{Quantifying Intersection} \label{s0203}

In the previous section, we have shown an infinite number of approximate intersections between a fractal and a line under certain conditions. This implies the problem of determining how much of the fractal falls within a translational or angular neighborhood of the line. Clearly what we need is a probability measure which can be easily evaluated in practice, so we utilize the well-known invariant measure in an explicit form.

\begin{dfn} \label{s020301}
We say that a measure $\mu:\R^d\rightarrow [0,1]$ is invariant with respect to the IFS $\mcT=\{T_1,\ldots,T_n\}$ if it satisfies the following equation
\[ \mu(S) = w_1\mu(T_1^{-1}(S)) + \ldots + w_n\mu(T_n^{-1}(S)) \]
for any $S\subset\R^d$ in its domain, with some fixed weights $w_k\in [0,1],\ \sum_k w_k=1$.
\end{dfn}

When the IFS maps are similitudes of the form $T_k(z):=p_k+\vl_kR_k(z-p_k)$ with $\vl_k\in (0,1)$ and unitary $R_k\in\R^{d\times d}$, then there exists a unique $s>0$, called the similarity dimension of the fractal, for which $\sum_k \vl_k^s=1$. In such a case, the weights $w_k=\vl_k^s$ seem natural.

\begin{thm} \label{s020302} (Hutchinson\cite{ba00007})
With respect to any IFS and weights, there exists a unique invariant probability measure with bounded support. We call this the invariant measure with respect to the IFS. For any initial probability measure with bounded support $\mu_0$, the recursion $\mu_{L}=w_1\mu_{L-1}\co T_1^{-1}+\ldots +w_n\mu_{L-1}\co T_n^{-1},\ L\in\N$ tends to the invariant measure as $L\tti$.
\end{thm}

\begin{thm} \label{s020303}
For any IFS $\mcT=\{T_1,\ldots,T_n\}$, with weights $w_1,\ldots,w_n$, and primary fixed point $p\in\mcP$, the invariant measure $\nu$ has the form
\[ \nu(S) = \lim_{L\tti} \sum(w_a : a\in\mcA,\ |a|=L,\ T_a(p)\in S)\ \ (S\subset\R^d)\ \ \mathrm{where}\ \ w_a:=\prod_{k=1}^{|a|} w_{a(k)}. \]
\end{thm}
\prf
The functions $\nu_0:=\delta_p$ and $\nu_L(S):=\sum(w_a : a\in\mcA,\ |a|=L,\ T_a(p)\in S),\ L\in\N, S\subset\R^d$ clearly satisfy the required properties of a measure, considering that
\[ 1 = (w_1+\ldots +w_n)^L = \sum(w_a : a\in\mcA,\ |a|=L)\ \ (L\in\N) \]
The Dirac measure $\nu_0$ has the support $\{p\}$ which is bounded. We now show that the weighted recursion in Theorem \ref{s020302} holds for $\nu_L,\ L\in\N$. The key to its derivation, is breaking up the sum according to the first coordinate of each address, which is taken from $\mcN=\{1,\dots,n\}$.
\[ \nu_L(S) = \sum(w_a : a\in\mcA,\ |a|=L,\ T_a(p)\in S) = \]
\[ = \sum_{k=1}^n \sum(w_a : a=(k,b)\in\mcN\times\mcA,\ |b|=L-1,\ T_k\co T_b(p)\in S) = \]
\[ = \sum_{k=1}^n w_k \sum(w_b : b\in\mcA,\ |b|=L-1,\ T_b(p)\in T_k^{-1}(S)) = \sum_{k=1}^n w_k\ \nu_{L-1}(T_k^{-1}(S)) \]
Note that the derivation also holds for $L=1$. By the above recursion and Theorem \ref{s020302} we have that $\exists\lim(\nu_L)=:\nu$ and it is the invariant measure with respect to the IFS $\mcT$. \sqr

The above explicit formula for the invariant measure has some practical advantages over other methods for its computation. One method well-known in image processing, is to discretize the plane over a rectangle, and approximate the invariant measure via a matrix recursion. Another method is carried out stochastically using Elton's Ergodic Theorem\cite{ba00031}. This exact formula for the invariant measure $\nu$ assigns weights $w_a$ to each fractal point $T_a(p)$ with an address of length $L$, and checks which points fall into the set $S$, then sums the weight of those points. Clearly $\nu$ can also be approximated in practice via its recursion up to $\nu_L$, with a large enough level $L\in\N$. We may thus find the density of intersections with the fractal in between parallel rays of light, reasonably spaced at $\ve:=\vl_{\min}^L\diam(\mcP)$, to see the intensity of the shadow in a particular direction, as illustrated on the figure below. This kind of analysis has clear applications to the design of fractal antennas for instance, and the tomography of fractalline structures.
\begin{figure}[H]
\centering
\includegraphics[width=470pt]{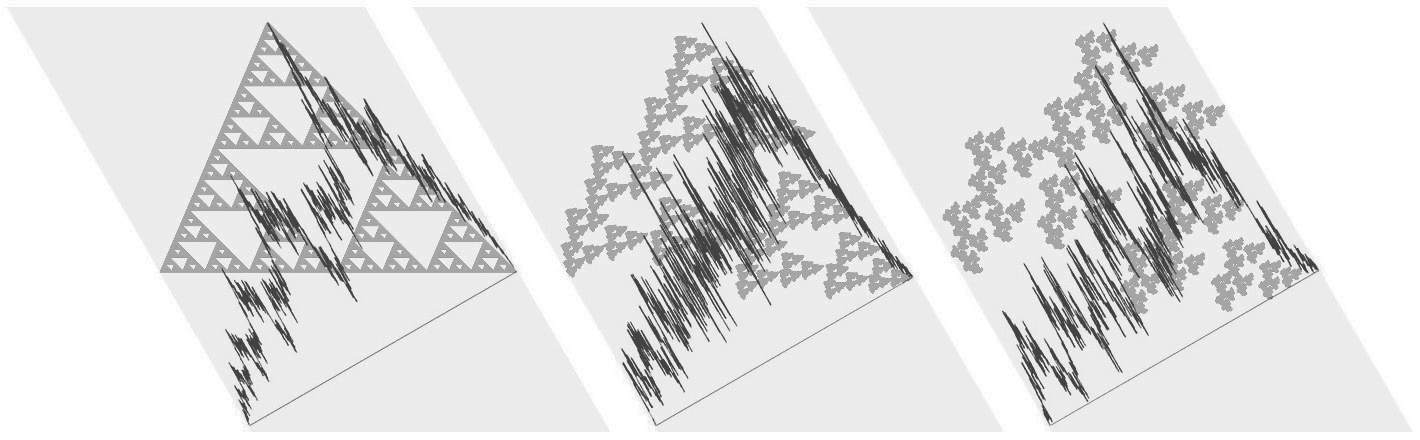}
\caption{Ray absorption density plots under rotational perturbation.}
\label{s020304}
\end{figure}

\section{Conclusions} \label{s03}

The paper highlighted the connection between the projection of IFS fractals and the Fractal-Line Intersection Problem, relevant to applications in computer graphics and antenna design. Broad classes of fractals have been introduced for which the projection from any direction is a segment. This was done in general, examining the intersection of hyperplanes and IFS fractals in $\R^d,\ d\in\N$. The method implies an algorithm for finding the points of intersection. The cardinality of intersection was discussed, and it was quantified via the invariant measure of IFS.

The author gratefully acknowledges the support of Prof. Edward R. Vrscay via NSERC.

\bibliographystyle{unsrt}
\bibliography{mybib}

\begin{thebibliography}{10}

\bibitem{bb00003}
Kenneth~J. Falconer.
\newblock {\em Fractal Geometry}.
\newblock Wiley, Chichester, second edition, 2003.

\bibitem{ba00018}
F.~Mendivil and T.~D. Taylor.
\newblock Thin sets with fat shadows: Projections of cantor sets.
\newblock {\em The American Mathematical Monthly}, 115:451--456, 2008.

\bibitem{bc00005}
John~C. Hart and Thomas~A. DeFanti.
\newblock Efficient antialiased rendering of 3-{D} linear fractals.
\newblock volume~25 of {\em Computer Graphics}, pages 91--100, Las Vegas, 1991.
  SIGGRAPH.

\bibitem{ba00030}
Nathan Cohen.
\newblock Fractal antennas part {I}.
\newblock {\em Communications Quarterly}, 9:7--22, 1995.

\bibitem{ba00019}
Robert~G. Hohlfeld and Nathan Cohen.
\newblock Self-similarity and the geometric requirements for frequency
  independence in antennae.
\newblock {\em Fractals}, 7:79--84, 1999.

\bibitem{bb00005}
Przemyslaw Prusinkiewicz and Aristid Lindenmayer.
\newblock {\em The Algorithmic Beauty of Plants}.
\newblock Springer-Verlag, second edition, 1996.

\bibitem{ba00020}
Boris Zeide.
\newblock Fractal analysis of foliage distribution in loblolly pine crowns.
\newblock {\em Canadian Journal of Forest Resources}, 28:106--114, 1998.

\bibitem{bb00001}
Benoit~B. Mandelbrot.
\newblock {\em The Fractal Geometry of {N}ature}.
\newblock Freeman, San Francisco, 1982.

\bibitem{ba00021}
A.~S. Besicovitch.
\newblock On the fundamental geometric properties of linearly measurable plane
  sets of points {III}.
\newblock {\em Math. Annalen}, 116:349--357, 1939.

\bibitem{ba00022}
H.~Federer.
\newblock The $(\varphi, k)$ rectifiable subsets of n-space.
\newblock {\em Transactions of the American Mathematical Society}, 62:114--192,
  1947.

\bibitem{ba00024}
J.~M. Marstrand.
\newblock Some fundamental geometrical properties of plane sets of fractional
  dimensions.
\newblock {\em Proc. Lond. Math. Soc.}, 4:257--302, 1954.

\bibitem{ba00025}
R.~Kaufman.
\newblock On the {H}ausdorff dimension of projections.
\newblock {\em Mathematika}, 15:153--155, 1968.

\bibitem{ba00026}
P.~Mattila.
\newblock Hausdorff dimension, orthogonal projections and intersections with
  planes.
\newblock {\em Ann. Acad. Sci. Fennicae}, A 1:227--244, 1975.

\bibitem{bb00006}
P.~Mattila.
\newblock {\em Geometry of Sets and Measures in Euclidean Spaces}.
\newblock Cambridge University Press, Cambridge, 1995.

\bibitem{ba00023}
R.~O. Davies.
\newblock On accessibility of plane sets and differentiation of functions of
  two real variables.
\newblock {\em Proc. Camb. Phil. Soc.}, 48:215--232, 1952.

\bibitem{ba00027}
Kenneth~J. Falconer.
\newblock Sets with prescribed projections and {N}ikodym sets.
\newblock {\em Proc. Lond. Math. Soc.}, 53:48--64, 1986.

\bibitem{ba00028}
K.~J. Falconer and J.~D. Howroyd.
\newblock Packing dimensions of projections and dimension profiles.
\newblock {\em Math. Proc. Cambridge Philos. Soc.}, 121:269--286, 1997.

\bibitem{ba00029}
J.~D. Howroyd.
\newblock Box and packing dimensions of projections and dimension profiles.
\newblock {\em Math. Proc. Cambridge Philos. Soc.}, 130:135--160, 2001.

\bibitem{ba00007}
John~E. Hutchinson.
\newblock Fractals and self similarity.
\newblock {\em Indiana University Mathematics Journal}, 30:713--747, 1981.

\bibitem{ba00017}
M.~F. Barnsley and S.~Demko.
\newblock Iterated function systems and the global construction of fractals.
\newblock {\em Proceedings of the Royal Society}, A399, 1985.

\bibitem{ba00031}
John~H. Elton.
\newblock An ergodic theorem for iterated maps.
\newblock {\em Ergodic Theory and Dynamical Systems}, 7:481--488, 1987.

\end{thebibliography}
\addcontentsline{toc}{chapter}{\textbf{References}}

\end{document}